\documentclass{birkmult} \usepackage{jms-bm}
\begin{document}
% Get revtex4 from http://publish.aps.org/revtex4

%\documentclass[amsmath,twocolumn,secnumarabic,superscriptaddress,%
%    floatfix,amssymb,aps,nofootinbib,showkeys]{revtex4}
%\documentclass[onesided]{birkmult}
%\usepackage{jms-bm}
%\usepackage{url,marginnote}
%\usepackage{graphicx,amsthm}
%\usepackage[percent]{overpic}
%\usepackage{varioref}
%\usepackage{placeins}
%\bibliographystyle{hamsalpha}

%\parskip0.5ex plus 0.1ex

%\newcommand{\arc}[1]{#1} %might want something like \widehat{#1}
\newcommand{\approaches}{\to}
\newcommand{\cross}{\times}
\newcommand{\union}{\cup}
\newcommand{\disj}{\sqcup} %not coprod, which is an upside down Pi
\newcommand{\bigdisj}{\bigsqcup}
\newcommand{\len}{\mathrm{Len}}
\newcommand{\mass}{\mathrm{Mass}}
\newcommand{\tc}{\mathrm{TC}}
\newcommand{\tce}{\mathrm{TC^*}}
\newcommand{\mesh}{\mathrm{Mesh}}
\newcommand{\eps}{\varepsilon}
\newcommand{\isom}{\cong}
\newcommand{\R}{\mathbb{R}}
\newcommand{\Z}{\mathbb{Z}}
\newcommand{\I}{\mathbb{I}}
\newcommand{\E}{\mathbb{E}}
\newcommand{\K}{\mathcal{K}}
\let\Section=\S
\renewcommand{\S}{\mathbb{S}}
\renewcommand{\d}{\partial}
\newcommand{\esssup}{\mathrm{ess\,sup}}
\renewcommand{\limsup}{\varlimsup}
\renewcommand{\liminf}{\varliminf}

\newcommand{\tild}[1]{\tilde{#1}}
\newcommand{\ttild}[1]{\hat{#1}}
\newcommand{\tg}{{\tild\gamma}}
\renewcommand{\phi}{\varphi}

\def\BV/{\text{\small{BV}}}
\def\FTC/{\text{\small{FTC}}}
\def\CAT/{\text{\small{CAT}}}

\renewcommand{\lem}[1]{Lemma~\ref{lem:#1}}
\newcommand{\thrm}[1]{Theorem~\ref{thm:#1}}
\renewcommand{\prop}[1]{Proposition~\ref{prop:#1}}
\renewcommand{\cor}[1]{Corollary~\ref{cor:#1}}
\newcommand{\xmpl}[1]{Example~\ref{ex:#1}}
\newcommand{\eqn}[1]{\eqref{eq:#1}}

% Other similar cross-ref commands:
\renewcommand{\figurename}{Figure}
\def\figr#1{Figure~\ref{fig:#1}}
\def\secn#1{Section~\ref{sec:#1}}

%%%%%%%% Including figures

%\graphicspath{{./Figs/}}
\newcommand{\incgr}[3]{\begin{overpic}[scale=#1]{#2}#3\end{overpic}}
\newcommand{\incgrw}[2]{\includegraphics[width=#1]{#2}}
\newcommand{\figcap}[3]{\caption[#2]{#3}\label{fig:#1}}

% For a single figure needing rescaling, use
%   \figs{label}{scale}{filename}{figure title}{caption}
% (where the title would only appear if we had a list of figures)
\newcommand{\figs}[6]{% 
  \begin{figure}
    \centerline{\incgr{#2}{#3}{#4}}
    \figcap{#1}{#5}{#6}
  \end{figure}
}
\newcommand{\figtwo}[4]{% 
  \begin{figure}
    \centerline{ \incgrw{.3\textwidth}{#2A} \hspace{.1\textwidth}
                 \incgrw{.3\textwidth}{#2B} }
    \figcap{#1}{#3}{#4}
  \end{figure}
}

%%%%%%%%%%%%%%

%Various options for our enumerate environment "enuma",
%starting with the one that seems best to me.  -jms

\renewcommand{\theenumi}{\alph{enumi}}
\renewcommand{\labelenumi}{(\theenumi)}
\newenvironment{enuma} {\begin{enumerate}}{\end{enumerate}}
\newcommand{\itm}[1]{({\ref{itm:#1}})}

%%%%%%%%%%%%%%%%%%%%%%%%%%%%ray, segment%%%%%%%%%%%%%%%%%%%%%%%%%%
\makeatletter
\newbox\overbox
\def\fakeover#1{\setbox\overbox\hbox{$#1$}\hbox
                         {$\overline{#1\hskip-\wd\overbox}$\hskip\wd\overbox}}
\def\overnoarrow#1{\mskip1.5mu\overline{\mskip-1.5mu#1}}
\def\overrightarrow#1{\mskip2mu\vbox{\m@th\ialign{##\crcr
   \rightarrowfill\crcr
   \noalign{\kern-.4pt               %compensate for thickness of rule
        \kern-\fontdimen22\textfont2 %pull down arrow so axis is on baseline
        \nointerlineskip}
   ${\mskip0mu\hfil\fakeover{#1}\hfil\mskip6mu}$\crcr}}\mskip-2mu}
\def\overleftrightarrow#1{\mskip-5mu\vbox{\m@th\ialign{##\crcr
   \leftarrowfill\hskip-.6em\rightarrowfill\crcr
   \noalign{\kern-.4pt               %compensate for thickness of rule
        \kern-\fontdimen22\textfont2 %pull down arrow so axis is on baseline
        \nointerlineskip}
   ${\mskip5mu\hfil\fakeover{#1}\hfil\mskip4mu}$\crcr}}\mskip-2mu}

%Note: the \smash makes TeX not leave extra space above the line.
%This probably wouldn't be a good idea without the \baselinestretch
%\def\baselinestretch{1.1}
\def\ray#1{{\smash{\overrightarrow{#1}}}}
\def\sline#1{{\smash{\overleftrightarrow{#1}}}}
\def\segment#1{{\smash{\overnoarrow{#1}}}}
\makeatother
%%%%%%%%%%%%%%%%%%%%%%%%%%%%%%%%%%%%%%%%%%%%%%%%%%%%%%%%%%%%%%%%%%%

% This is some code to loosen up how LaTeX places floats:
% helpful for a heavily illustrated paper.

%\renewcommand{\topfraction}{.85}
%\renewcommand{\bottomfraction}{.7}
%\renewcommand{\textfraction}{.15}
%\renewcommand{\floatpagefraction}{.66}
%\renewcommand{\dbltopfraction}{.66}
%\renewcommand{\dblfloatpagefraction}{.66}
%\setcounter{topnumber}{9}
%\setcounter{bottomnumber}{9}
%\setcounter{totalnumber}{20}
%\setcounter{dbltopnumber}{9}
%
%%\begin{document}

\title{Curves of Finite Total Curvature}

\author{John M. Sullivan}
\address{\tuaddr32}
\email{sullivan@math.tu-berlin.de}


\begin{abstract}
We consider the class of curves of finite total curvature,
as introduced by Milnor.  This is a natural class for
variational problems and geometric knot theory, and
since it includes both smooth and polygonal curves,
its study shows us connections between discrete and
differential geometry.  To explore these ideas, we consider
theorems of F{\'a}ry/Milnor, Schur, Chakerian and Wienholtz.
%These notes are from a lecture delivered at the 
%Oberwolfach Seminar \emph{Discrete Differential Geometry}
%held the first week of June 2004.
\end{abstract}

\maketitle

\bigskip\noindent
Here we introduce the ideas of discrete differential geometry
in the simplest possible setting: the geometry and curvature of curves,
and the way these notions relate for polygonal and smooth curves.
\index{curvature (of curve)}%
The viewpoint has been partly inspired by work in geometric knot theory, 
which studies geometric properties of space curves in relation to
their knot type, and looks for optimal shapes for given knots.

After reviewing Jordan's definition of the length of a curve,
we consider Milnor's analogous definition~\cite{Milnor} of total curvature.
In this unified treatment, polygonal and smooth curves are
both contained in the larger class of \FTC/ (finite total curvature) curves.
%\index{curve!FTC|seeonly{finite total curvature}}%
\index{FTC|seeonly{finite total curvature}}\index{finite total curvature (FTC)!curve}%
We explore the connection between \FTC/ curves and \BV/ functions.
Then we examine the theorems of F{\'a}ry/Milnor, Schur and Chakerian in terms
of \FTC/ curves.  We consider relations between total curvature and
Gromov's distortion,
and then we sketch a proof of a result by Wienholtz in integral geometry.
We end by looking at ways to define curvature density for polygonal curves.

A companion article~\cite{DS-ConvFTC} examines more carefully
the topology of \FTC/ curves, showing that any two
sufficiently nearby \FTC/ graphs are isotopic.
The article~\cite{Sul-CSDS}, also in this volume, looks at
curvatures of smooth and discrete surfaces; the discretizations
are chosen to preserve various integral curvature relations.

Our whole approach in this survey should be compared to that
of Alexandrov and Reshetnyak~\cite{AlexResh}, who develop
much of their theory for curves having one-sided tangents everywhere,
a class somewhat more general than \FTC/.

\section{Length and total variation}

We want to consider the geometry of curves.  Of course curves---unlike
higher-di\-men\-sion\-al manifolds---have no local \emph{intrinsic} geometry.
So we mean the \emph{extrinsic} geometry of how the curve sits in some
ambient space~$M$.  Usually $M$ will be in euclidean $d$-space~$\E^d$,
but the study of space curves naturally leads also to the study of curves on
spheres.  Thus we also allow~$M$ to be a smooth Riemannian manifold;
for convenience we embed~$M$ isometrically into some~$\E^d$.
(Some of our initial results would still hold with~$M$ being any path-metric
space; compare~\cite{AlexResh}.  Here, however, our curves will be quite
arbitrary but not our ambient space.)

A curve is a one-dimensional object, so we start by recalling the
topological classification of one-manifolds:
A compact one-manifold (allowing boundary) is a finite
disjoint union of components, each homeomorphic to an interval $\I:=[0,L]$
or to a circle $\S^1:=\R/L\Z$.
Then a \emph{parametrized curve} %\index{curve!parametrized}
in~$M$ is a continuous map from a compact one-manifold to~$M$.
That is, each of its components is a (parametrized) \emph{arc}
$\gamma:\I\to M$ or \emph{loop} $\gamma:\S^1\to M$.
A loop can be viewed as an arc whose endpoints are equal and thus identified.
A \emph{curve} in~$M$ is an equivalence class of parametrized
curves, where the equivalence relation is given by orientation-preserving
reparametrization of the domain.

(An unoriented curve would allow arbitrary reparametrization.
Although we will not usually care about the orientation of our curves,
keeping it around in the background is convenient, fixing for
instance a direction for the unit tangent vector of a rectifiable curve.)

Sometimes we want to allow reparametrizations by arbitrary
monotone functions that are not necessarily homeomorphisms.
Intuitively, we can collapse any time interval on which
the curve is constant, or conversely stop for some time
at any point along the curve.  Since there might be
infinitely many such intervals, the easiest formalization
of these ideas is in terms of Fr{\'e}chet distance~\cite{Frechet}.

The \emph{Fr{\'e}chet distance}\index{Frechet distance@Fr{\'e}chet distance}
between two curves is the infimum, over all stricly monotonic
reparametrizations, of the maximum pointwise distance.
(This has also picturesquely been termed, perhaps originally
in~\cite{AltGod}, the ``dog-leash distance'':
the minimum length of leash required for a dog who walks forwards along one
curve while the owner follows the other curve.)
Two curves whose Fr{\'e}chet distance is zero are
equivalent in the sense we intend: homeomorphic
reparametrizations that approach the infimal value zero will limit
to the more general reparametrization that might eliminate or
introduce intervals of constancy.  See also \cite[\Section~X.7]{Graves}.

Given a connected parametrized curve~$\gamma$, a choice of
$$0\le t_1 < t_2 < \cdots < t_n \le L$$
gives us the vertices $v_i:=\gamma(t_i)$ of an \dfn{inscribed polygon}~$P$,
whose edges are the minimizing geodesics $e_i:=v_iv_{i+1}$
in~$M$ between consecutive vertices.
(If $\gamma$ is a loop, then indices~$i$ are to be taken modulo~$n$,
that is, we consider an inscribed polygonal loop.)
We will write $P<\gamma$ to denote that $P$ is a polygon inscribed in~$\gamma$.

The edges are uniquely determined by the vertices when $M=\E^d$
or more generally whenever~$M$ is simply connected with
nonpositive sectional curvature (and thus a $\CAT/(0)$ space).
When minimizing geodesics are not unique, however,
as in the case when~$M$ is a sphere and two consecutive vertices
are antipodal, some edges may need to be separately specified.
The \emph{mesh}\index{mesh (of polygon)} of~$P$
(relative to the given parametrization of~$\gamma$) is
$$\mesh(P) := \max_i(t_{i+1}-t_i).$$
(Here, of course, for a loop the $n^{\text{th}}$ value
in this maximum is $(t_1+L)-t_n$.)

The \emph{length} of a polygon is simply the sum of the edge lengths:
$$\len(P)=\len_M(P):=\sum_i d_M(v_i,v_{i+1}).$$
This depends only on the vertices and not on which
minimizing geodesics have been picked as the edges, since by
definition all minimizing geodesics have the same length.
If $M\subset N\subset \E^d$, then a given set of
vertices defines (in general) different polygonal curves
in~$M$ and~$N$, with perhaps greater length in~$M$.

We are now ready to define the \emph{length}\index{length (definition)}
of an arbitrary curve:
$$\len(\gamma) := \sup_{P<\gamma} \len(P).$$
When $\gamma$ itself is a polygonal curve, it is easy to check that
this definition does agree with the earlier one for polygons.
This fact is essentially the definition of what it means
for~$M$ to be a path metric space: the distance $d(v,w)$ between any
two points is the minimum length of paths connecting them.
By this definition, the length of a curve $\gamma\subset M\subset\E^d$
is independent of~$M$; length can be measured in~$\E^d$ since even
though the inscribed polygons may be different in~$M$, their supremal
length is the same.

This definition of length for curves originates
with Jordan~\cite{Jordan} and independently Scheeffer~\cite{Scheeffer},
and is thus often called ``\ix{Jordan length}''.
(See also \cite[\Section 2]{Cesari}.)
For $C^1$-smooth curves it
can easily be seen to agree with the standard integral formula.

\begin{lemma}
Given a polygon~$P$, if~$P'$ is obtained
from~$P$ by deleting one vertex~$v_k$
then $\len(P')\le\len(P)$.
We have equality here if and only if
$v_k$ lies on a minimizing geodesic from~$v_{k-1}$ to~$v_{k+1}$.
\end{lemma}

\begin{proof}
This is simply the triangle inequality applied to
the triple $v_{k-1}$, $v_k$, $v_{k+1}$.
\end{proof}

A curve is called \emph{rectifiable}\index{rectifiable curve}
if its length is finite.  (From the beginning, we have considered only
compact curves.  Thus we do not need to distinguish rectifiable and
locally rectifiable curves.)

\begin{proposition}
A curve is rectifiable if and only if it admits a Lipschitz parametrization.
\end{proposition}
\begin{proof}
If $\gamma$ is $K$-Lipschitz on $[0,L]$, then its length is at most $KL$,
since the Lipschitz bound gives this directly for any inscribed polygon.
Conversely, a rectifiable curve can be reparametrized by its arclength
$$s(t):=\len\big(\gamma|_{[0,t]}\big)$$
and this arclength parametrization is $1$-Lipschitz.
\end{proof}

If the original curve was constant on some time
interval, the reparametrization here will not be one-to-one.
A nonrectifiable curve has no Lipschitz parametrization,
but might have a H{\"o}lder-continuous one.
(For a nice choice of parametrization for an arbitrary curve,
see~\cite{Morse-param}.)

Given this proposition, the standard theory of Lipschitz functions shows
that a rectifiable curve~$\gamma$ has almost everywhere a well-defined
unit tangent vector $T=\gamma'$, its derivative with respect
to its arclength parameter~$s$.  Given a rectifiable curve,
we will most often use this arclength parametrization. 
The domain is then $s\in[0,L]$ or $s\in\R/L\Z$, where $L$ is the length.

Consider now an arbitrary function~$f$ from $\I$ (or $\S^1$) to~$M$,
\emph{not} required to be continuous.  We can apply the same
definition of inscribed polygon~$P$, with vertices $v_i=f(t_i)$,
and thus the same definition of length $\len(f)=\sup\len(P)$.
This length of~$f$ is usually called the \dfn{total variation} of~$f$,
and~$f$ is said to
be~\emph{\BV/} (of \emph{bounded variation})\index{bounded variation (BV)}
when this is finite.

For a fixed ambient space $M\subset\E^d$,
the total variation of a discontinuous~$f$ as a function to~$M$
may be greater than its total variation in~$\E^d$.
The supremal ratio here is
$$\sup_{p,q\in M} \frac{d_M(p,q)}{d_{\E^d}(p,q)},$$
what Gromov called the \dfn{distortion} of the
embedding $M\subset\E^d$.  (See \cite[pp.~6--9]{GLP},
\cite[p.~114]{GromovFRM} and \cite{GromovHED}, as well
as \cite{KS-disto, DS-disto}.)
When $M$ is compact and smoothly embedded (like~$\S^{d-1}$),
this distortion is finite; thus $f$ is \BV/ in~$M$ if and
only if it is \BV/ in~$\E^d$.

The class of \BV/ functions (here, from $\I$ to~$M$) is often useful
for variational problems.  Basic facts about \BV/ functions can be
found in the original book~\cite{Carath} by Cara\-th{\'e}o\-dory
or in many analysis texts like
%\cite[Sect.~5.2]{Royden},
\cite[Sect.~2.19]{GoffPed}, \cite[Chap.~3]{Boas} or~\cite{Berberian}.
For more details and higher dimensions,
see for instance \cite{Ziemer, AFP}.

Here, we recall one nice characterization:
$f$ is \BV/ if and only if it has a weak (distributional) derivative.
Here, a weak derivative means an $\E^d$-valued Radon measure~$\mu$\index{measure}
which plays the role of $f'\,dt$ in integration by parts, meaning that
$$\int_0^L f \phi'\,dt = - \int_0^L \phi\,\mu$$
for every smooth test function~$\phi$ vanishing at the endpoints.
(This characterization of \BV/ functions is
one form of the Riesz representation theorem.)

\begin{proposition}
If $f$ is \BV/, then $f$ has well-defined right and left limits
$$f_\pm(t) := \lim_{\tau\to t^{\pm}} f(\tau)$$
everywhere.  Except at countably many \emph{jump points} of~$f$,
we have $f_-(t)=f(t)=f_+(t)$.
\end{proposition}
\begin{proof}[Sketch of proof]
We consider separately each of the $d$ real-valued coordinate
functions~$f^i$.
We decompose the total variation of~$f^i$ into positive and negative
parts, each of which is bounded.  This lets us write $f^i$ as the
difference of two monotonically increasing functions. (This is its
so-called \emph{Jordan decomposition}.)  An increasing
function can only have countably many (jump) discontinuities.
(Alternatively, one can start by noting that a real-valued
function without, say, a left-limit at~$t$ has infinite total
variation even locally.)
\end{proof}

In functional analysis, \BV/ functions are often
viewed as equivalence classes of functions differing only on sets
of measure zero.  Then we replace total variation with \emph{essential
total variation}, the infimal total variation over the equivalence
class.  A minimizing representative will be necessarily continuous
wherever $f_-=f_+$.  A unique representative can be obtained by
additionally requiring left (or right, or upper, or lower) semicontinuity
at the remaining jump points.

Our definition of curve length is not very practical,
being given in terms of a supremum over all possible inscribed polygons.
But it is easy to find a sequence of polygons
guaranteed to capture the supremal length:
\begin{proposition}\label{prop:lenlimit}
Suppose $P_k$ is a sequence of polygons inscribed in a curve~$\gamma$ such that
$\mesh(P_k)\to0$.  Then $\len(\gamma)=\lim\len(P_k)$.
\end{proposition}
\begin{proof}
By definition, $\len(\gamma)\ge\len(P_k)$, so $\len(\gamma)\ge\limsup(P_k)$.
Suppose that $\len(\gamma)>\liminf(P_k)$.  Passing to a subsequence,
for some~$\eps>0$ we have $\len(\gamma)\ge\len(P_k)+2\eps$.
Then by the definition of length, there is an inscribed~$P_0$
(with, say, $n$ vertices) such that $\len(P_0)\ge\len(P_k)+\eps$
for all~$k$.

The common refinement of~$P_0$ and~$P_k$ is of course
at least as long as~$P_0$.
But this refinement is $P_k$ with a fixed number~($n$)
of vertices inserted;
for each~$k$, these $n$ insertions together add length at least~$\eps$
to~$P_k$.  For large enough~$k$, these $n$ insertions are
at disjoint places along~$P_k$, so their effect on the length
is independent of the order in which they are performed.
Passing again to a subsequence, there is thus
some vertex $v_0=\gamma(t_0)$ of~$P_0$ such that if~$P^0_k$
is~$P_k$ with~$v_0$ inserted, we have $\len(P^0_k)\ge\len(P_k)+\eps/n$.

But $\gamma$ is continuous, in particular at~$t_0$.  So there exists
some $\delta>0$ such that for $t\in[t_0-\delta,t_0+\delta]$ we have
$d_M(\gamma(t),\gamma(t_0))<\eps/2n$.  Choosing $k$ large enough
that $\mesh(P_k)<\delta$, the vertices of~$P_k$ immediately
before and after~$v_0$ will be within this range,
so $\len(P^0_k)<\len(P_k)+\eps/n$, a contradiction.
\end{proof}

Although we have stated this proposition only for continuous curves~$\gamma$,
the same holds for \BV/ functions~$f$, as long as~$f$ is semicontinuous
at each of its jump points.

An analogous statement does not hold for polyhedral approximations to
surfaces.  First, an inscribed polyhedron (whose vertices lie ``in order''
on the surface) can have greater area than the original surface, even if
the mesh size (the diameter of the largest triangle) is small.
Second, not even the limiting value is guaranteed to be correct.
Although Serret had proposed~\cite[p.~293]{Serret} defining
surface area as a limit of polyhedral areas, claiming this limit
existed for smooth surfaces, Schwarz soon found a counterexample,
now known as the ``\ix{Schwarz lantern}''~\cite{Schwarz-lantern}:
seemingly nice triangular meshes inscribed in a cylinder, with mesh size
decreasing to zero, can have area approaching infinity.

Lebesgue~\cite{Lebesgue} thus defined surface area
as the \emph{lim$\,$inf} of such converging polyhedral areas.
(See \cite{AlpTor,Cesari-Chern} for an extensive discussion of
related notions.)
One can also rescue the situation with the additional requirement
that the shapes of the triangles stay bounded (so that their
normals approach that of the smooth surface), but we will not
explore this here.  (See also~\cite{Tonelli}.
In this volume, the companion article~\cite{Sul-CSDS} treats
curvatures of smooth and discrete surfaces, and~\cite{xWardetzky}
considers convergence issues.)

Historically, such difficulties led to new approaches to
defining length and area, such as Hausdorff measure.
These measure-theoretic approaches work well in all dimensions,
and lead to generalizations of submanifolds like the currents
and varifolds of geometric measure theory (see~\cite{morgan}).
We have chosen here to present the more ``old-fashioned'' notion of
Jordan length for curves because it nicely parallels Milnor's definition
of total curvature, which we consider next.

\section{Total curvature}\label{sec:tc}
Milnor~\cite{Milnor} defined a notion of total curvature for
arbitrary curves in euclidean space.
Suppose $P$ is a polygon in~$M$ with no two consecutive vertices equal.
Its \dfn{turning angle} at an interior vertex~$v_n$
is the angle $\phi\in[0,\pi]$ between the oriented tangent vectors at~$v_n$
to the two edges~$v_{n-1}v_n$ and~$v_nv_{n+1}$.
(Here, by saying interior vertices, we mean to exclude the
endpoints of a polygonal arc, where there is no turning angle; 
every vertex of a polygonal loop is interior.
The supplement of the turning angle, sometimes called
an interior angle of~$P$, will not be of interest to us.)

If~$M$ is an oriented surface, for instance if $M=\E^2$ or~$\S^2$,
then we can also define a \emph{signed turning angle}
$\phi\in[-\pi,\pi]$ at~$v_n$,
except that where $\phi=\pm\pi$ its sign is ambiguous.

To find the \dfn{total curvature} $\tc(P)$ of a polygon~$P$,
we first collapse any sequence of consecutive equal vertices to a single vertex.
Then $\tc(P)$ is simply the sum of the turning angles at all interior vertices.

Here, we mainly care about the case when~$P$ is in $M=\E^d$.
Then the unit tangent vectors along the edges,
in the directions $v_{n+1}-v_n$,
are the vertices of a polygon in~$\S^{d-1}$ called the
\dfn{tantrix} of~$P$.  (The word is a shortening of ``tangent indicatrix''.)
The total curvature of~$P$ is the length of its tantrix in~$\S^{d-1}$.

\begin{lemma}
\emph{(See \cite[Lemma 1.1]{Milnor} and~\cite{Borsuk}.)}
Suppose~$P$ is a polygon in~$\E^d$.
If~$P'$ is obtained from~$P$ by deleting one vertex~$v_n$
then $\tc(P')\le\tc(P)$.  We have equality here if
$v_{n-1}v_nv_{n+1}$ are collinear in that order, or if
$v_{n-2}v_{n-1}v_nv_{n+1}v_{n+2}$ lie convexly in some two-plane,
but never otherwise.
\end{lemma}
\begin{proof}
Deleting~$v_n$ has the following effect on the tantrix:
two consecutive vertices (the tangents to the edges
$v_{n-1}v_n$ and $v_nv_{n+1}$) are coalesced into a single one
(the tangent to the edge $v_{n-1}v_{n+1}$).
It lies on a great circle arc connecting the original two,
as in \figr{tc-change}.
Using the triangle inequality twice,
the length of the tantrix decreases (strictly, unless
the tantrix vertices $v_{n-1}v_n$ and $v_nv_{n+1}$ coincide, or
the relevant part lies along a single great circle in~$\S^{d-1}$).
\end{proof}

\figs{tc-change}{.5}{tc-change}
{\put(31,85){$v_{n-2}v_{n-1}$}
 \put(49,69){$v_{n-1}v_n$}
 \put(37,51){$v_{n-1}v_{n+1}$}
 \put(10,45){$v_n v_{n+1}$}
 \put(50,17){$v_{n+1}v_{n+2}$} }
{Deleting vertex decreases total curvature}
{Four consecutive edges of a polygon~$P$ in space give four vertices
and three connecting edges (shown here as solid lines) of
its tantrix on the sphere.  When the middle vertex~$v_n$ of~$P$ is deleted,
two vertices of the tantrix get collapsed to
a single new one (labeled $v_{n-1}v_{n+1}$); it lies
somewhere along edge connecting the two original vertices.
The two new edges of the tantrix are shown as dashed lines.
Applying the triangle inequality twice, we see that the length
of the new tantrix (the total curvature of the new polygon) is no greater.}

\begin{corollary}\label{cor:polytc}
If~$P$ is a polygon in~$\E^d$ and $P'<P$ then $\tc(P')\le\tc(P)$.
\end{corollary}

\begin{proof}
Starting with~$P$, first insert the vertices of~$P'$;
since each of these lies along an edge of~$P$,
these insertions have no effect on the total curvature.
Next delete the vertices not in~$P'$;
this can only decrease the total curvature.
\end{proof}

\begin{defn}
For any curve $\gamma\subset\E^d$ we follow Milnor~\cite{Milnor} to define
$$\tc(\gamma) := \sup_{P<\gamma} \tc(P).$$
\index{finite total curvature (FTC)!curve}%
We say that $\gamma$ has \emph{finite total curvature}
(or that $\gamma$ is \emph{\FTC/}) if $\tc(\gamma)<\infty$.
\end{defn}
When $\gamma$ is itself a polygon,
this definition agrees with the first one by \cor{polytc}.
Our curves are compact, and thus lie in bounded subsets of $\E^d$.
It is intuitively clear then that a compact curve of infinite length must
have infinite total curvature; that is, that all \FTC/ curves are rectifiable.
This follows rigorously by applying the quantitative estimate
of \prop{len-bd} below to finely inscribed polygons,
using Propositions~\ref{prop:lenlimit} and~\ref{prop:tclimit}.

Various properties follow very easily from this definition.
For instance, if the total curvature of an arc is less than~$\pi$ then
the arc cannot stray too far from its endpoints.   In particular,
define the \emph{spindle} of angle~$\theta$ with endpoints~$p$
and~$q$ to be the body of revolution bounded by a circular arc of total
curvature~$\theta$ from~$p$ to~$q$ that has been revolved about~$\sline{pq}$.
(The spindle is convex, looking like an American football, for $\theta\le\pi$
and is a round ball for $\theta=\pi$.)

\begin{lemma}
Suppose $\gamma$ is an arc from~$p$ to~$q$ of total curvature
$\phi<\pi$.  Then $\gamma$ is contained in the spindle of
angle~$2\phi$ from~$p$ to~$q$.
\end{lemma}

\begin{proof}
Suppose $x\in\gamma$ is outside the spindle.  Consider the
planar polygonal arc~$pxq$ inscribed in~$\gamma$.
In that plane, since $x$ is outside the circular arc of total curvature~$2\phi$,
by elementary geometry the turning angle of~$pxq$ at~$x$ is
greater than~$\phi$, contradicting the definition of~$\tc(\gamma)$.
\end{proof}

We also immediately recover \ix{Fenchel's theorem}~\cite{Fenchel}:

\begin{theorem}[Fenchel]\label{thm:fenchel}
Any closed curve in~$\E^d$ has total curvature at least~$2\pi$.
\end{theorem}
\begin{proof}
Pick any two distinct points $p$,~$q$ on the curve.
(We didn't intend the theorem to apply to the constant curve!)
The inscribed polygonal loop from~$p$ to~$q$ and back
has total curvature~$2\pi$, so the original curve has
at least this much curvature.
\end{proof}

For this approach to Fenchel's theorem to be satisfactory,
we do need to verify (as Milnor did~\cite{Milnor})
that our definition of total curvature agrees with the usual
one $\int\!\kappa\,ds$ for smooth curves.  For us, this will
follow from \prop{tclimit}.

\section{First variation of length}
We can characterize \FTC/ curves as those with \BV/ tangent vectors.
This relates to the variational characterization of curvature
in terms of first variation of length.
(The discussion in this section is based on \cite[Sect.~4]{CFKSW1}.)

We have noted that an \FTC/ curve~$\gamma$ is rectifiable,
hence has a tangent vector~$T$ defined almost everywhere.
We now claim that the total curvature of~$\gamma$ is exactly the length
(or, more precisely, the essential total variation)
of this tantrix~$T$ as a curve in~$\S^{d-1}$.
We have already noted this for polygons, so the general case
seems almost obvious from the definitions.
However, while the tantrix of a polygon inscribed in~$\gamma$ is a spherical
polygon, it is not inscribed in~$T$; instead its vertices are
averages of small pieces of~$T$.  Luckily, this is close enough
to allow the argument of \prop{lenlimit} to go through again:
Just as for length, in order to compute total curvature it suffices
to take any limit of finer and finer inscribed polygons.

\begin{proposition}\label{prop:tclimit}
Suppose~$\gamma$ is a curve in~$\E^d$.
If $P_k$ is a sequence of polygons inscribed in~$\gamma$
with $\mesh(P_k)\to0$, then $\tc(\gamma)=\lim\tc(P_k)$.
This equals the essential total variation
of its tantrix $T\subset\S^{d-1}$.
\end{proposition}
We leave the proof of this proposition as an exercise.
The first statement essentially follows as in the
proof of \prop{lenlimit}: if it failed there would be
one vertex~$v_0$ along~$\gamma$ whose insertion would cause
a uniform increase in total curvature for all polygons
in a convergent subsequence, contradicting the fact
that sufficiently small arcs before and after~$v$
have arbitrarily small total curvature.  The second statement
follows by measuring both $\tc(\gamma)$ and the total
variation through limits of (different but nearby) fine polygons.

To summarize, a rectifiable curve~$\gamma$ has finite total curvature
if and only if its unit tangent vector $T=\gamma'(s)$
is a function of bounded variation.  (Thus the space of \FTC/
curves could be called  $W^{1,BV}$ or $BV^1$.)
If $\gamma$ is \FTC/, it follows that at every point of~$\gamma$ there are
well-defined left and right tangent vectors $T_{\pm}$; these are equal
and opposite except at countably many points, the \emph{corners} of~$\gamma$.

Now, to investigate curvature from a variational point of view,
suppose we consider a \emph{continuous deformation} $\gamma_t$ of a
curve $\gamma\subset\E^d$: fixing any parametrization of~$\gamma$, this means
a continuous family $\gamma_t$ of parametrized curves with $\gamma_0=\gamma$.
(If we reparametrize~$\gamma$, we must apply the same reparametrization
to each~$\gamma_t$.)

We assume that position of each point is (at least~$C^1$) smooth in time;
the initial velocity of~$\gamma_t$ will then be given by
some (continuous, $\E^d$-valued) \emph{vectorfield} $\xi$ along~$\gamma$.

Let $\gamma$ be a rectifiable curve parametrized by arclength~$s$,
with unit tangent vector~$T=\gamma'(s)$ (defined almost everywhere).
Suppose $\gamma_t$ is a variation of $\gamma=\gamma_0$ whose
initial velocity~$\xi(s)$ is a Lipschitz function of arclength.
Then the arclength derivative
$\xi'=\d\xi/\d s$ is defined almost everywhere along~$\gamma$,
and a standard first-variation calculation shows that
$$\delta_\xi \len(\gamma) := \frac{d}{dt}\Big|_0{\len(\gamma_t)}
     = \int_\gamma \big<T,\xi'\big>\,ds.$$
If $\gamma$ is smooth enough, we can integrate this by parts to get
$$\delta_\xi \len(\gamma)
= -\int_\gamma \big<T',\xi\big>\,ds
  - \sum_{x\in\d\gamma} \big<\pm T,\xi\big>,$$
where, in the boundary term, the sign is chosen
to make $\pm T$ point inward at~$x$.
In fact, not much smoothness is required: 
as long as $\gamma$ is \FTC/, we know that its unit
tangent vector $T$ is \BV/, so we can interpret $T'$ as a
measure, and this first-variation
formula holds in the following sense:
the weak (distributional) derivative $\K:=T'$ 
is an $\E^d$-valued Radon measure
along~$\gamma$ which we call the \dfn{curvature force}.

On a $C^2$ arc of~$\gamma$,  the curvature force is
$\K = d T = \kappa N\,ds$ and is absolutely continuous
with respect to the arclength or Hausdorff measure $ds=\mathcal{H}^1$.
The curvature force has an atom (a point mass or Dirac delta)
at each corner $x\in \gamma$, with $\K\{x\}=T_+(x)+T_-(x)$.
Note that at such a corner, the mass of~$\K$ is not
the turning angle~$\theta$ at~$x$.  Instead,
$$|\K|\{x\} = |\K\{x\}| = 2\sin(\theta/2).$$

Therefore, the total mass (or total variation) $|\K|(\gamma)$
of the curvature force~$\K$ is somewhat less than the
total curvature of~$\gamma$: at each corner it counts
$2\sin(\theta/2)$ instead of~$\theta$.  Whereas $\tc(\gamma)$
was the length (or total variation) $\len_{\S^{d-1}}(T)$ of
the tantrix~$T$ viewed as a (discontinuous) curve on the sphere~$\S^{d-1}$,
we recognize this total mass as its length $\len_{\E^d}(T)$
in euclidean space.  Thus we call it the euclidean total curvature
of~$\gamma$:
$$\tce(\gamma) := \len_{\E^d}(T) = |\K|(\gamma).$$

%When $\gamma$ has boundary, we choose to include in $\K$ an atom
%at each endpoint of $\gamma$, with mass~$1$ and
%pointing in the inward tangent direction.  This means we
%need no boundary terms in the formula
%$\delta_\xi \len(\gamma) = -\int_\gamma \big<\xi,\K\big>$.

Returning to the first variation of length,
we say that a vectorfield $\xi$ along $\gamma$
is \emph{smooth} if $\xi(s)$ is a smooth function of arclength.
The first variation $\delta \len(\gamma)$ can be viewed as a linear
functional on the space of smooth vectorfields $\xi$ along~$\gamma$.
As such a \emph{distribution}, it has degree zero, by definition, if 
$\delta_\xi \len(\gamma)=\int_\gamma \big<T,\xi'\big>\,ds$ is bounded by
$C\sup_\gamma|\xi|$ for some constant $C$.  This happens exactly when we
can perform the integration by parts above.

\pagebreak
We collect the results of this section as:
\begin{proposition}\label{prop:ftc}
Given any rectifiable curve $\gamma$, the following conditions are equivalent:
\begin{enuma}
\item $\gamma$ is $\FTC/$.
\item There exists a curvature force $\K=dT$ along $\gamma$ such that
  $$\delta_\xi \len(\gamma) =
    - \int_\gamma \big<\xi,\K\big> - \sum_{\d\gamma} \big<\xi,\pm T\big>.$$
\item The first variation $\delta \len(\gamma)$ has distributional degree zero.
\qed
\end{enuma}
\end{proposition}

Of course, just as not all continuous functions are \BV/,
not all $C^1$ curves are \FTC/.  However, given
an \FTC/ curve, it is piecewise~$C^1$ exactly when it has finitely
many corners, and is $C^1$ when it has no corners,
that is, when $\K$ has no atoms.
The \FTC/ curve is furthermore $C^{1,1}$ when $T$ is Lipschitz,
or equivalently when $\K$ is absolutely continuous (with respect
to arclength~$s$) and has bounded Radon/Nikodym derivative
$d\K/ds = \kappa N$.

\section{Total curvature and projection}
The F\'ary/Milnor theorem\index{Fary/Milnor@F\'ary/Milnor theorem}
says that a knotted curve in~$\E^3$
has total curvature at least~$4\pi$, twice that of
an unknotted round circle.  The different proofs given
by F\'ary~\cite{fary} and Milnor~\cite{Milnor} can both be
interpreted in terms of a proposition about the average total curvature
of different projections of a curve.\index{projection!of curve}

The \ix{Grassmannian} $G_k\E^d$ of $k$-planes in $d$-space
is compact, with a unique rotation-invariant probability measure~$d\mu$.
For $p\in G_k\E^d$, we denote by~$\pi_p$ the orthogonal projection
to~$p$.  When we speak about averaging over all projections,
we mean using~$d\mu$.  This proposition is essentially due
to F{\'a}ry~\cite{fary}, though he only stated the case $d=3$, $k=2$.

\begin{proposition}\label{prop:tc-ave}
Given a curve~$K$ in~$\E^d$, and some fixed $k<d$,
the total curvature of~$K$ equals the average total curvature
of its projections to $k$-planes.  That is,
$$\tc(K)=\int_{G_k\E^d} \!\!\!\! \tc\big(\pi_p(K)\big) \,d\mu.$$
\end{proposition}

\begin{proof}
By definition of total curvature and \prop{tclimit},
we may reduce to the case where $K$ is a polygon.
(To interchange the limit of ever finer inscribed polygons
with the average over the Grassmannian, we use
the Lebesgue monotone convergence theorem.)
Since the total curvature of a polygon is the sum of
its turning angles, it suffices to consider a single
angle.  So let~$P_\theta$ be a three-vertex polygonal
arc with a single turning angle of~$\theta\in[0,\pi]$.
Defining
$$f_k^d(\theta):=\int_{p\in G_k\E^d} \!\!\!\! \tc\big(\pi_p(P_\theta)\big)\,d\mu,$$
the rotation-invariance of $\mu$ shows this is independent
of the position of $P_\theta$, and our goal is to show
$f_k^d(\theta)=\theta$.

First note that $f_k^d$ is continuous.  It is also additive:
$$f_k^d(\alpha+\beta)= f_k^d(\alpha)+f_k^d(\beta)$$
when $0\le\alpha,\beta\le\alpha+\beta\le\pi$.
This follows by cutting the single corner of~$P_\theta$
into two corners of turning angles~$\alpha$ and~$\beta$.
Any projection of the resulting convex planar arc is still convex and planar,
so additivity holds in each projection, and thus also holds after averaging.

A continuous, additive function is linear: $f_k^d(\theta)=c_k^d\theta$
for some constant $c_k^d$.  Thus we just need to show $c_k^d=1$.
But we can easily evaluate $f$ at $\theta=\pi$,
where $P_\pi$ is a ``cusp'' in which the incoming
and outgoing edges overlap.  Clearly every projection of $P_\pi$
is again such a cusp with turning angle~$\pi$
(except for a set of measure zero where the projection is a single point).
Thus $f_k^d(\pi)=\pi$, so $c_k^d=1$ and we are done.
\end{proof}

A curve in~$\E^1$ has only cusps for corners, so its
total curvature will be a nonnegative multiple of~$\pi$.  A loop
in~$\E^1$ has total curvature a positive multiple of~$2\pi$.
(In particular, the loop in~$\E^1$ is a real-valued function on~$\S^1$,
and its total curvature is $2\pi$ times the number of local maxima.)

This proposition could be used to immediately reduce the
$d$-dimensional case of Fenchel's theorem (\thrm{fenchel})
to the $k$-dimensional case. 
Historically, this could have been useful.  For instance, the theorem
is trivial in~$\E^1$ by the previous paragraph.  In~$\E^2$, the idea
that a simple closed curve has total signed curvature $\pm2\pi$ essentially
dates back to Riemann (compare~\cite[\Section 1]{Chern}).
Fenchel's proof~\cite{Fenchel} (which was his 1928
doctoral dissertation in Berlin) was for~$d=3$, and the first
proof for general dimensions seems to be that of Borsuk~\cite{Borsuk}.

For alternative proofs of Fenchel's theorem---as well as comparisons
among them---see~\cite{Horn-Fenchel} and~\cite[\Section 4]{Chern},
and the references therein,
especially~\cite{Liebmann-Fenchel} and~\cite{MR12:634d}. %Fenchel, BAMS
Voss~\cite{Voss-Fenchel} related the total curvature of a space curve
to the total Gauss curvature of a tube around it, and thereby gave
new a new proof of the F{\'a}ry/Milnor theorem as well as Fenchel's
theorem.

Milnor's proof~\cite{Milnor} of the F\'ary/Milnor theorem can be rephrased as
a combination of the case $k=1$ of \prop{tc-ave}
with the following:

\begin{lemma}
If~$K\subset\E^3$ is nontrivially knotted, then
any projection of~$K$ to~$\E^1$ has total curvature at least~$4\pi$.
\end{lemma}
\begin{proof}[Sketch of proof]
If there were some projection direction in which the total
curvature were only~$2\pi$, then the corresponding (linear) height
function on~$\E^3$ would have only one local minmum and one
local maximum along~$K$.  Then at each intermediate height,
there are exactly two points of~$K$.  Connecting each such pair
with a straight segment, we form a spanning disk
showing $K$ is unknotted.
\end{proof}

(This proof isn't quite complete as stated: at some intermediate
heights, one or even both strands of~$K$ might have a whole
subarc at that constant height.  One could still patch in a disk,
but easier is to follow Milnor and at the very beginning
replace $K$ by an isotopic inscribed polygon.  Compare~\cite{DS-ConvFTC}.)

F\'ary~\cite{fary}, on the other hand, used $k=2$ in \prop{tc-ave}, having proved:

\begin{lemma}
Any nontrivial knot diagram (a projection to~$\E^2$
of a knotted curve $K\subset\E^3$)
has total curvature at least~$4\pi$.
\end{lemma}

\begin{proof}[Sketch of proof]
Any knot diagram divides the plane into regions: one unbounded
and several bounded.  If every bounded region is adjacent to
the unbounded region, the only possibility is a tree-like
diagram as in \figr{fary}~(left); this is clearly unknotted no matter
how we choose over- and under-crossings.

\figtwo{fary}{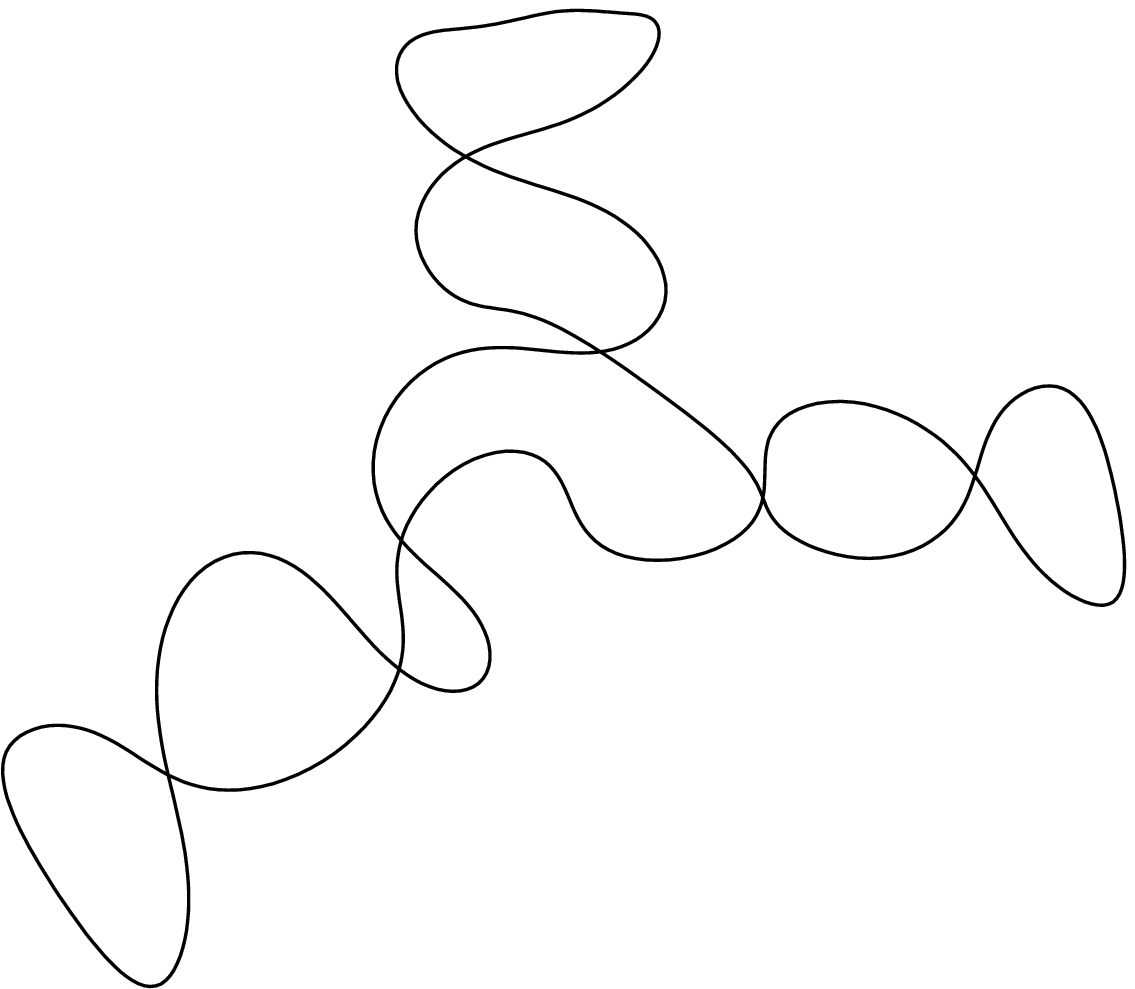}{F{\'a}ry's proof}
{A diagram in which every region is adjacent to the outside
is in fact unknotted (left), so a knot diagram has a doubly
enclosed region (right).}

Thus every nontrivial knot diagram~$D$ has some region $R$ which
is ``doubly enclosed'' by the curve, not necessarily in the
sense of oriented winding numbers, but in the sense that any
ray outwards from a point $p\in R$ must cut the curve twice.
Then $R$ is part of the \emph{second hull} of the curve~\cite{CKKS},
and the result follows by Lemmas~5 and~6 there. To summarize the
arguments there (which parallel those of~\cite{fary}), note first that 
the cone over~$D$ from~$p$ has cone angle at least~$4\pi$ at~$p$;
by Gauss/Bonnet, this cone angle equals the total signed geodesic
curvature of~$D$ in the cone, which is at most its total (unsigned) curvature.
\end{proof}

Either of these lemmas, combined with the appropriate case
of \prop{tc-ave}, immediately yields the
F\'ary/Milnor theorem\index{Fary/Milnor@F\'ary/Milnor theorem}.

\begin{corollary}
The total curvature of a nontrival knot $K\subset\E^3$ is at least~$4\pi$.
\qed
\end{corollary}

While F\'ary stated the appropriate case of \prop{tc-ave}
pretty much as such, we note that Milnor didn't speak about
total curvature of projections to lines, but instead only
about extrema of height functions.  (A reinterpretation more
like ours can already be found, for instance, in~\cite{Taniyama-tcgraph}.)

Denne~\cite{Denne} has given a beautiful new proof of F{\'a}ry/Milnor
analogous to the easy proof we gave for Fenchel's \thrm{fenchel}.
Indeed, the hope that there could be such a proof had led us to conjecture
in~\cite{CKKS} that every knot~$K$ has an \emph{alternating
quadrisecant}.  A quadrisecant, by definition, is a line in space intersecting
the knot in four points~$p_i$.  The~$p_i$ form an inscribed polygonal loop,
whose total curvature (since it lies in a line) is either~$2\pi$ or~$4\pi$,
depending on the relative ordering of the~$p_i$ along the line
and along~$K$.  The quadrisecant is called alternating exactly when the
curvature is~$4\pi$.  Denne proved this conjecture: every nontrivially
knotted curve in space has an alternating quadrisecant.
This gives as an immediate corollary not only the F{\'a}ry/Milnor theorem,
but also a new proof that a knot has a \emph{second hull}~\cite{CKKS}.

We should also note that a proof of F{\'a}ry/Milnor
for knots in hyperbolic as well as euclidean space
was given by Brickell and Hsiung~\cite{BriHsi}; essentially
they construct a point on the knot that also lies in its second hull.
The theorem is also now known for knots in an arbitrary Hadamard manifold,
that is, a simply connected manifold of nonpositive curvature:
Alexander and Bishop~\cite{AlxBsh:FM} find a finite sequence of
inscriptions---first a polygon~$P_1$ inscribed in~$K$, then
inductively~$P_{i+1}$ inscribed in~$P_i$---ending with a $pqpq$ quadrilateral,
while Schmitz~\cite{Schmitz} comes close to constructing a quadrisecant.

The results of this section would not
be valid for $\tce$ in place of $\tc$: there is
no analog to \prop{tc-ave}, and the F\'ary/Milnor
theorem would fail.

We conclude this section by recalling that another standard proof
of Fenchel's theorem uses the following integral-geometric lemma due
to Crofton (see \cite[\Section 4]{Chern} and also \cite{Santalo-Chern, Santalo})
to conclude that a spherical curve of length less than~$2\pi$
is contained in an open hemisphere:
\begin{lemma} \label{lem:sphcrof}
The length of a curve $\gamma\subset\S^{d-1}$ equals $\pi$ times
the average number of intersections of $\gamma$ with great
hyperspheres~$\S^{d-2}$.
\end{lemma} 
\begin{proof}
It suffices to prove this for polygons (appealing again
to the monotone convergence thereom).
Hence it suffices to consider a
single great-circle arc.  But clearly for such arcs,
the average number of intersections is proportional
to length.  When the length is~$\pi$, (almost) every
great circle is intersected exactly once.
\end{proof}

This lemma, applied to the tantrix, is equivalent
to the case $k=1$ of our \prop{tc-ave}.  Indeed,
when projecting~$\gamma$ to the line in direction~$v$,
the total curvature we see counts the number of
times the tantrix intersects the great sphere
normal to~$v$.  We note also that knowing this case $k=1$
(for all~$d$) immediately implies all other cases of \prop{tc-ave},
since a projection from~$\E^d$ to~$\E^1$ can be
factored as projections $\E^d\to\E^k\to\E^1$.

Finally, we recall an analogous statement of the
famous \index{Cauchy Crofton@Cauchy/Crofton formula}{Cauchy/Crofton formula}.
Its basic idea dates back to Buffon's 1777 analysis~\cite{Buffon}
of his needle problem (compare \cite[Chap.~1]{KlainRota}).
Cauchy obtained the formula by 1841~\cite{Cauchy-area} and generalized
it to find the surface area of a convex body.
Crofton's 1868 paper~\cite{Crofton} on geometric probability
includes this among many integral-geometric formulas for plane curves.
(See~\cite{AyaDub} for a treatment like ours
for rectifiable curves,
and also \cite[\Section 1.7C]{doCarmo} and~\cite{Santalo-Chern}.)
\begin{lemma}[Cauchy/Crofton]\label{lem:crofton}
The length of a plane curve equals $\pi/2$ times
the average length of its projections to lines.
\end{lemma}
\begin{proof}
Again, we can reduce first to polygons, then to
a single segment.  So the result certainly holds
for some constant; to check the constant is~$\pi/2$,
it is easiest to compute it for the unit circle, where
every projection has length~$4$.
\end{proof}

In all three of our integral-geometric arguments
(\ref{prop:tc-ave}, \ref{lem:sphcrof}, \ref{lem:crofton})
we proved a certain function was linear and then found
the constant of proportionality by computing one
(perhaps sometimes surprising) example.  One finds
also in the literature proofs where the integrals (over the circle
or more generally the Grassmannian) are simply computed explicitly.
Although the trigonometric integrals are not too difficult,
that approach seems to obscure the geometric essence of the argument.

\section{Schur's comparison theorem}
\ix{Schur's comparison theorem}~\cite{Schur} is a well-known result saying
that straightening an arc will increase the distance between its
endpoints.
 
Chern, in \Section 5 of his beautiful essay~\cite{Chern} in his 
MAA book, gives a proof for $C^2$ curves and remarks (without
proof) that it also applies to piecewise smooth curves.
In~\cite{CKS2} we noted that Chern's proof actually applies
to $C^{1,1}$ curves, that is to curves with a Lip\-schitz tangent
vector, or with bounded curvature density.  In fact, the
natural class of curves to which the proof applies is \FTC/ curves.

\begin{theorem}[Schur's Comparison Theorem]\label{thm:schur}
Let $\bar\gamma\subset\E^2$ be a planar arc such that
joining the endpoints of $\bar\gamma$ results in a convex
(simple, closed) curve, and let $\gamma\subset\E^d$ be
an arc of the same length~$L$.
Suppose that~$\bar\gamma$ has nowhere less curvature than~$\gamma$
with respect to the common arclength parameter $s\in[0,L]$,
that is, that for any subinterval $I\subset[0,L]$ we have
$$\tc\big(\bar\gamma|_I\big) \ge \tc\big(\gamma|_I\big).$$
(Equivalently, $|\K|-|\bar\K|$ is a nonnegative measure.)
Then the distance between the endpoints is greater for~$\gamma$:
$$\big|\gamma(L)-\gamma(0)\big| \ge \big|\bar\gamma(L)-\bar\gamma(0)\big|.$$
\end{theorem}

\begin{proof}
By convexity, we can find an~$s_0$ such that
the (or, in the case of a corner, some supporting)
tangent direction~$T_0$ to~$\bar\gamma$ at~$\bar\gamma(s_0)$
is parallel to $\bar\gamma(L)-\bar\gamma(0)$.
Note that the convexity assumption implies that the total
curvature of either half of~$\bar\gamma$ (before or after~$s_0$)
is at most~$\pi$.
Now move $\gamma$ by a rigid motion so that it shares this same
tangent vector at~$s_0$, as in \figr{schur}.
\figs{schur}{1.25}{schur}
{\put(50,48){$\gamma(s_0)$}
 \put(79,48){$T_0$}
 \put(9,3){$\bar\gamma(0)$}
 \put(69,3){$\bar\gamma(L)$}
 \put(-5,28){$\gamma(0)$}
 \put(87,9){$\gamma(L)$}
}
{Proof of Schur's Theorem}
{Schur's Theorem compares the space curve~$\gamma$ (solid)
with a planar curve~$\bar\gamma$ (dotted) that has pointwise no less
curvature.  The proof places them with a common tangent vector~$T_0$
(at arclength~$s_0$)
pointing in the direction connecting the endpoints of~$\bar\gamma$.}
(If $\gamma$ has a corner at~$s_0$, so does $\bar\gamma$.
In this case we want to arrange that the angle from $T_0$
to each one-sided tangent vector $T_\pm$ is
at least as big for~$\bar\gamma$ as for~$\gamma$.)

By choice of~$T_0$ we have
$$\big|\bar\gamma(L)-\bar\gamma(0)\big|
= \big<\bar\gamma(L)-\bar\gamma(0),T_0\big>
= \int_0^L \!\big<\bar T(s),T_0\big>\,ds,$$
while for~$\gamma$ we have
$$\big|\gamma(L)-\gamma(0)\big|
\ge \big<\gamma(L)-\gamma(0),T_0\big>
= \int_0^L \!\big<T(s),T_0\big>\,ds.$$
Thus it suffices to prove (for almost every~$s$) that
$$\big<T(s),T_0\big> \ge \big<\bar T(s),T_0\big>.$$

But starting from~$s_0$ and moving outwards in either direction,
$\bar T$ moves straight along a great circle arc, at speed
given by the pointwise curvature; in total it moves less than
distance~$\pi$.  At the same time, $T$ moves
at the same or lower speed, and perhaps not straight but on
a curved path.  Clearly then $T(s)$ is always closer to~$T_0$
than $\bar T(s)$ is, as desired.  More precisely,
\begin{equation*}
\big<\bar T(s),T_0\big>= \cos\tc\big(\bar\gamma|_{[s_0,s]}\big)
\le \cos\tc\big(\gamma|_{[s_0,s]}\big) \le \big<T(s),T_0\big>.\qedhere
\end{equation*}
\end{proof}

The special case of Schur's theorem when $\gamma$
and $\bar\gamma$ are polygons is usually
called \index{Cauchy arm lemma}{Cauchy's arm lemma}.
It was used in Cauchy's proof~\cite{Cauchy} of the rigidity
of convex polyhedra, although Cauchy's own proof of the arm lemma
was not quite correct, as discovered 120 years later by Steinitz.
The standard modern proof of the arm lemma
(due to Schoenberg; see~\cite{AigZieg} or~\cite[p.~235]{Cromwell-polyh})
is quite different from the proof we have given here.
For more discussion of the relation between Schur's theorem and Cauchy's lemma,
see \cite{Connelly-invent,ORourke-Schur}.

The history of this result is somewhat complicated.
Schur~\cite{Schur} considered only
the case where $\gamma$ and~$\bar\gamma$ have pointwise
equal curvature: twisting a convex plane curve out of the plane
by adding torsion will increase its chord lengths.  He
considers both polygonal and smooth curves.  He attributes
the orginal idea (only for the case where~$\bar\gamma$ is
a circular arc) to unpublished work of H.~A.~Schwarz in 1884.
The full result, allowing the space curve to have less
curvature, is evidently due to Schmidt~\cite{Schmidt-Schur}.
See also the surveys by Blaschke in~\cite{Blaschke1}
and \cite[\Section 28--30]{Blaschke2}.

In Schur's theorem, it is irrelevant whether we
use the spherical or euclidean version of total curvature.
If we replace $\tc$ by $\tce$ throughout, the statement and
proof remain unchanged, since the curvature comparison is pointwise.

\ixsection{Chakerian's packing theorem}

A less familiar result due to Chakerian (and cited for instance
as \cite[Lemma~1.1]{BuckSimon}) captures the intuition
that a long rope packed into a small ball must have large curvature.
\begin{proposition}\label{prop:len-bd}
A connected \FTC/ curve~$\gamma$ contained in the unit
ball in~$\E^d$ has length no more than $2+\tce(\gamma)$.
(If~$\gamma$ is closed, the~$2$ can be omitted.) 
\end{proposition}
\begin{proof}
Use the arclength parametrization $\gamma(s)$.  Then
\begin{align*}
\len(\gamma) &= \int 1\,ds = \int \big<T,T\big>\,ds \\
 &= \big<T,\gamma\big>\big|_{\text{endpts}} - \int\big<\gamma,d\K\big> \\
 &\le 2 + \int d|\K| = 2+\tce(\gamma).\qedhere
\end{align*}
\end{proof}

Chakerian~\cite{Chak2} gave exactly this argument
for $C^2$ curves and then used a limit argument (rounding the
corners of inscribed polygons) to get a version for all curves.
Note, however, that this limiting procedure gives the
bound with~$\tce$ replaced by~$\tc$;
this is of course equivalent for $C^1$ curves but weaker for
curves with corners.  For closed curves, Chakerian noted that
equality holds in $\len\le\tc$ only for a great circle
(perhaps traced multiple times).
In our sharper bound $\len\le\tce$, we have equality
also for a regular $n$-gon inscribed in a great circle.

(We recall that we appealed to this theorem in \secn{tc} to
deduce that \FTC/ curves are rectifiable.  This is not circular reasoning:
we first apply the proof above to polygons, then deduce that
\FTC/ curves are rectifiable and indeed have \BV/ tangents,
and fianlly apply the proof above in general.)

Chakerian had earlier~\cite{Chak1} given a quite different
proof (following F{\'a}ry) that $\len\le\tc$.
We close by interpreting that first argument in our framework.
Start by observing that in~$\E^1$, where curvature is quantized,
it is obvious that for a closed curve in the unit ball
(which is just a segment of length~$2$)
$$\len\le\tce=\tfrac2{\pi}\,\tc.$$
Combining this with Cauchy/Crofton (\lem{crofton}) and
our \prop{tc-ave} gives immediately $\len\le\tc$ for curves in
the unit disk in~$\E^2$.  With a little care, the same is true
for curves that fail to close by some angular holonomy.
(The two endpoints are at equal radius, and we do include
in the total curvature the angle they make when they are
rotated to meet.)  Rephrased, the length of a
curve~$\gamma$ in a unit neighborhood of the cone point on
a cone surface of arbitrary cone angle is at most the
total (unsigned) geodesic curvature of~$\gamma$ in the cone.
Finally, given any curve~$\gamma$ in the unit ball in~$\E^d$,
Chakerian considers the cone over~$\gamma$ from the origin.
The length is at most the total curvature in the cone,
which is at most the total curvature in space.  Rather
than trying to consider cones over arbitrary \FTC/ curves,
we can prove the theorem for polygons and then take a limit.

\section{Distortion}

We have already mentioned Gromov's distortion for an embedded
submanifold. For a curve~$\gamma\subset\E^d$, the \dfn{distortion} is
$$\delta(\gamma):=\sup_{p\ne q\in\gamma} \delta(p,q), \qquad
\delta(p,q) := \frac{\len(p,q)}{|p-q|},$$
where $\len(p,q)$ is the (shorter) arclength distance along~$\gamma$.
Here, we discuss some relations between distortion and total
curvature; many of these appeared in the first version of~\cite{DS-disto},
but later improvements to the main argument there made the discussion
of \FTC/ curves unnecessary.

Examples like a steep logarithmic spiral show that
arcs of infinite total curvature can have finite distortion,
even distortion arbitrarily close to~$1$.   However, there is an easy
bound the other way:

\begin{proposition}
Any arc of total curvature $\alpha<\pi$ has distortion
at most $\sec(\alpha/2)$.
\end{proposition}

\begin{proof}
First, note that it suffices to prove this for the endpoints of the arc.
(If the distortion were realized by some other pair $(p,q)$, we
would just replace the original arc by the subarc from~$p$ to~$q$.)

Second, note by that Schur's \thrm{schur} we may assume the arc is
convex and planar: we replace any given arc by the
locally convex planar arc with the same pointwise curvature.
Because the total curvature is less than~$\pi$, the planar
arc is globally convex in the sense of \thrm{schur}, and
the theorem shows the endpoint separation has only decreased.

Now fix points $p$ and~$q$ in the plane; for any given tangent
direction at~$p$, there is a unique triangle~$pxq$ with
exterior angle~$\alpha$ at~$x$.
Any convex arc of total curvature~$\alpha$
from~$a$ (with the given tangent) to~$c$
lies within this triangle.
By the Cauchy/Crofton formula of \lem{crofton},
its length is then at most that of the polygonal arc~$pxq$.
Varying now the tangent at~$a$, the locus of points~$x$ is a circle,
and it is easy to see that the length is maximized in the symmetric
situation, with $\delta=\sec(\alpha/2)$.
\end{proof}

This result might be compared with the bound \cite[Lemma~5.1]{KS-disto}
on distortion for a $C^{1,1}$ arc with bounded curvature
density~$\kappa\le1$.  By Schur's \thrm{schur} such an arc
of length $2a\le2\pi$ can be compared to a circle, and
thus has distortion at most $a/\sin a$.

For any curve~$\gamma$, the distortion is realized either by a pair
of distinct points or in a limit as the points approach,
simply because $\gamma\times\gamma$ is compact.
In general, the latter case might be quite complicated.
On an \FTC/ curve, however, we now show that the distortion between
nearby pairs behaves very nicely.
Define $\alpha(r)\le\pi$ to be the turning angle
at the point $r\in\gamma$, with $\alpha=0$ when $r$ is not a corner.

\begin{lemma}\label{lem:disto-sc}
On an \FTC/ curve~$\gamma$, we have
$$\limsup_{p,q\to r} \delta(p,q) = \sec\frac{\alpha(r)}2,$$
with this limit realized by symmetric pairs $(p,q)$ approaching~$r$
from opposite sides.
\end{lemma}
\begin{proof}
The existence of one-sided tangent vectors~$T_{\pm}$ at~$r$
is exactly enough to make this work, since the quotient in
the definition of $\delta(p,q)$ is similar to the difference
quotients defining~$T_{\pm}$.  Indeed,
near~$r$ the curve looks very much like a pair of rays
with turning angle~$\alpha$.  Thus the lim$\,$sup is the same
as the distortion of these rays, which is $\sec(\alpha/2)$,
realized by any pair of points symmetrically spaced about the vertex.
\end{proof}

This leads us to define $\delta(r,r):=\sec(\alpha(r)/2)$, giving
a function $\delta:\gamma\times\gamma\to[1,\infty]$ that is
upper semicontinuous.  The compactness of $\gamma$
then immediately gives:
\begin{corollary}\label{cor:disto-real}
On an FTC curve~$\gamma$, there is a pair $(p,q)$ of
(not necessarily distinct) points on~$\gamma$ which realize
the distortion $\delta(\gamma)=\delta(p,q)$.
\qed\end{corollary}

Although distortion is not a continuous functional
on the space of rectifiable curves, it is lower semicontinuous.
A version of the next lemma appeared as~\cite[Lem.~2.2]{KS-disto}:
\begin{lemma} \label{lem:c0-lim}
Suppose curves $\gamma_j$ approach a limit $\gamma$ in the sense
of Fr{\'e}chet distance.  Then 
$\delta(\gamma) \le \liminf \delta(\gamma_j)$.
\end{lemma}
\begin{proof}
The distortion for any fixed pair of points is lower semicontinuous
because the arclength between them is.  (And length will indeed
jump down in a limit unless the tangent vectors also converge
in a certain sense.  See \cite[Chap.~2, \Section 29]{Tonelli}.)
The supremum of a family of lower semicontinuous functions
is again lower semicontinuous.
\end{proof}

\section{A projection theorem of Wienholtz}
\index{Wienholtz's projection theorem}

In~\cite{KS-disto}, we made the following conjecture:

\emph{
Any closed curve~$\gamma$ in~$\E^d$ of length~$L$ has some
orthogonal projection to~$\E^{d-1}$ of diameter at most~$L/\pi$.
}

This yields an easy new proof of Gromov's result
(see \cite{KS-disto,DS-disto}) that a closed curve
has distortion at least~$\pi/2$, that of a circle.
Indeed, consider the height function along~$\gamma$
in the direction on some projection of small diameter.
For any point $p\in\gamma$, consider the \emph{antipodal}
point~$p^*$ halfway around~$\gamma$, at arclength~$L/2$.
Since the height difference between~$p$ and~$p^*$ is
continuous and changes sign, it equals zero for some $(p,p^*)$.
The distance between the projected images of these points is at most
the diameter, at most~$L/\pi$.  But since the heights were equal,
this distance is the same as their distance $|p-p^*|$ in~$\E^d$.

For $d=2$, we noted that our conjecture follows immediately
from Cauchy/Crofton: a closed plane curve of length~$L$
has average width~$L/\pi$ and thus has
width at most~$L/\pi$ in some direction.
But for higher~$d$, the analogs of Cauchy/Crofton
give a weaker result.  (A curve of length~$L$
in~$\E^3$, for instance, has projections to~$\E^2$
of average length $\pi L/4$, and thus has some
planar projection of diameter at most~$\pi L/8$.)

In a series of Bonn preprints from~1999, Daniel Wienholtz
proved our conjecture and in fact somewhat more:
a closed curve in~$\E^d$ of length~$L$ has some orthogonal
projection to~$\E^{d-1}$ which lies in a ball of diameter~$L/\pi$.
Because Wienholtz's work has unfortunately remained unpublished,
we outline his arguments here.

\begin{proposition}\label{prop:antip}
Given any closed curve~$\gamma$ in~$\E^d$ for $d\ge3$,
there is some slab containing~$\gamma$, bounded by parallel
hyperplanes~$h_1$ and~$h_2$, with points $a_i,b_i\in\gamma\cap h_i$
occuring along~$\gamma$ in the order $a_1a_2b_1b_2$.
(We call the $h_i$ a pair of
parallel interleaved bitangent support planes for~$\gamma$.)
\end{proposition}

\begin{proof}
Suppose not.  Then for any unit vector $v\in\S^{d-1}$, we can divide
the circle parameterizing $\gamma$ into two complementary arcs
$\alpha(v)$, $\beta(v)$, such that the (global) maximum of the height function
in direction $v$ is achieved only (strictly) within $\alpha$,
and the minimum is achieved only (strictly) within $\beta$.
In fact, these arcs can be chosen to depend continuously on $v$.
Now let $a(v)\in\S^1$ be the midpoint of $\alpha(v)$.  Consider
the continuous map $a\colon \S^{d-1}\rightarrow \S^1 \subset \E^{d-1}$.
By one version of the Borsuk/Ulam theorem (see~\cite{Matousek}),
there must be some $v$ such that $a(v)=a(-v)$.  But the height functions
for $v$ and $-v$ are negatives of each other, so maxima in direction $v$
live in $\alpha(v)$ and in $\beta(-v)$, while minima live in $\beta(v)$ and
$\alpha(-v)$.  This is impossible if $a(v)$ is the midpoint of both $\alpha(v)$
and $\alpha(-v)$.  (In fact, a sensible choice for~$\alpha$ makes~$a$
antipodally equivariant: $a(-v)=-a(v)$, allowing direct application
of another version of Borsuk/Ulam.)
\end{proof}

\begin{lemma}\label{lem:projab}
If $\gamma$ is a curve in $\E^{m+n}$ of length $L$, and
its projections to $\E^m$ and $\E^n$ have lengths $a$ and $b$,
then $a^2 + b^2 \le L^2$.
\end{lemma}

\begin{proof}
By \prop{lenlimit}, it suffices to prove this for polygons.
Let $a_i,b_i\ge0$ be the lengths of the two projections of
the $i^{\textrm{th}}$ edge, and consider the polygonal arc in~$\E^2$
with successive edge vectors $(a_i,b_i)$.
Its total length is $\sum\sqrt{a_i^2+b_i^2}=L$,
but the distance between its endpoints is $\sqrt{a^2+b^2}$.
\end{proof}

\begin{theorem}
Any closed curve~$\gamma$ in~$\E^d$ of length~$L$
lies in a cylinder of diameter $L/\pi$.
\end{theorem}

\begin{proof}
As we have noted,
the case $d=2$ follows directly from the Cauchy/Crofton formula (\lem{crofton}),
since~$\gamma$ has width at most $L/\pi$ in some direction.
We prove the general case by induction.  So given $\gamma$ in~$\E^d$,
find two parallel interleaved bitangent support planes with normal~$v$,
as in \prop{antip}.  Let~$\tau$ be the distance between these planes,
the thickness of the slab in which $\gamma$ lies.
Project $\gamma$ to a curve $\bar\gamma$ in the plane orthogonal to $v$,
and call its length $\bar L$.  By induction, $\bar\gamma$ lies in a
cylinder of radius $\bar L/2\pi$.
Clearly, $\gamma$ lies in a parallel cylinder of radius $r$, centered in the
middle of the slab, where $r^2 = (\bar L/2\pi)^2 + (\tau/2)^2$.  So we need to
show that $(\bar L/2\pi)^2 + (\tau/2)^2 \le (L/2\pi)^2$, i.e.,
$\bar L^2 + \pi^2 \tau^2 \le L^2$.
In fact, since the length of the projection of~$\gamma$ to the
one-dimensional space in direction~$v$ is at least~$4\tau$,
\lem{projab} gives us $\bar L^2 + 4^2 \tau^2 \le L^2$,
which is better than we needed.
\end{proof}

If we are willing to settle for a slightly worse bound in the original
conjecture, Wienholtz also shows that we can project in a known direction:

\begin{proposition}\label{prop:wien2}
Suppose a closed curve $\gamma\subset\E^d$ has length~$L$,
and $p_1, p_2 \in \gamma$ are points realizing its diameter.
Then its projection $\bar\gamma$ to the plane orthogonal to $p_1-p_2$
has diameter at most $L/2\sqrt2$.  This estimate is sharp for a square.
\end{proposition}

\begin{proof}
Let $a_1$, $a_2$ be the preimages of a pair of points realizing
the diameter~$D$ of the projected curve~$\bar\gamma$.  We may
assume $\gamma$ is a quadrilateral with vertices
$a_1$, $p_1$, $a_2$, $p_2$, since any other curve would be longer.
(This reduces the problem to some
affine $\E^3$ containing these four points.)

Along~$\gamma$, first suppose the $a_i$ are interleaved with the $p_i$,
so that the quadrilateral is $a_1p_1a_2p_2$, as in \figr{wienholtz}~(left).
We can now reduce to $\E^2$:
rotating $a_1$, $a_2$ independently about the line~$\sline{p_1p_2}$
fixes the length, but maximizes the diameter $D$ when
the points are coplanar, with $a_i$ on opposite sides of the line.
Now let $R$ be the reflection across this line~$\sline{p_1p_2}$,
and consider the vector $R(p_1-a_1)+(a_1-p_2)+(p_1-a_2)+R(a_2-p_2)$.
Its length is at most~$L$,
but its component in the direction $p_1-p_2$
is at least~$2D$, as is its component in the perpendicular direction.
Therefore $L \ge 2\sqrt2 D$.

\figs{wienholtz}{.6}{wienholtz}
{\put(-2,27){$a_1$}
 \put(45.5,27){$a_2$}
 \put(22,51){$p_1$}
 \put(22,4){$p_2$}
 \put(75,51){$p_1$}
 \put(75,4){$p_2$}
 \put(79,31){$\theta$}
 \put(52,22){$a_1$}
 \put(98,32){$a_2$}
}
{The projection of a quadrilateral}
{These quadrilaterals show the sharp bounds
for the two cases in \prop{wien2}.
The original diameter is $p_1p_2$ in both cases,
and the projected diameter is~$D$ (shown at the bottom).
The square $a_1p_1a_2p_2$ (left) has length $2\sqrt2 D$.
The bowtie $a_1a_2p_1p_2$ (right) has length
$2D/\sin\theta + D/\cos\frac{\theta}2$, minimized at about $3.33$
for $\theta\approx76^\circ$ as shown, %76.3454
but in any case clearly greater than~$3D$.
}

Otherwise, write the quadrilateral as $a_1a_2p_1p_2$
as in \figr{wienholtz}~(right), and let $\theta$ be the angle between the
vectors $a_2-a_1$ and $p_1-p_2$.
Suppose (after rescaling) that $1=|a_2-a_1|\le|p_1-p_2|$.  Then
the projected diameter is $D=\sin\theta$.  By the triangle inequality,
the two remaining sides have lengths summing to at least
$$|p_1-a_2+a_1-p_2|\ge2\sin(\theta/2).$$
(Equality here holds for instance
when the quadrilateral is a symmetric bowtie.) Thus
\begin{equation*}
L/D\ge \frac{2+2\sin(\theta/2)}{\sin\theta}
      = \frac{2}{\sin\theta} + \frac1{\cos(\theta/2)}\ge 2+1
      > 2\sqrt2.\qedhere
\end{equation*}
\end{proof}

\section{Curvature density}

We have found a setting which treats polygons and smooth curves
in a unified way as two special cases of the more general
class of finite total curvature curves.  Many standard results
on curvature, like Schur's comparison theorem, work nicely
in this class.

However, there is some ambiguity in how to measure curvature
at a corner, reflected in our quantities~$\tc$ and~$\tce$.
A corner of turning angle~$\theta$ is counted either as~$\theta$
or as $2\sin(\theta/2)$, respectively.

At first, $\tc$ seems more natural: if we round off a corner
into any convex planar arc, its curvature is~$\theta$.
And the nice behavior of~$\tc$ under projection (\prop{tc-ave})
explains why it is the right quantity for results
like the F\'ary/Milnor theorem.

But from a variational point of view, $\tce$, which measures the mass
of the curvature force $\K$ at a corner, is sometimes most natural.
\prop{len-bd} is an example of a result whose
sharp form involves $\tce$.
An arbitrary rounding of a corner, whether or not it is convex or
planar, will have the same value of
$$\K=T_+-T_-=2\sin(\theta/2)\,N.$$
When we do choose a smooth, convex, planar rounding,
we note that
$$\int\big|\kappa N\big|\,ds=\int\!\kappa\,ds=\theta,
\qquad \text{while} \quad
\bigg|\int\!\kappa N\,ds\bigg| = 2\sin(\theta/2) .$$

For us, the curvature of an \FTC/ curve has been given
by the measure~$\K$.  For a polygon, of course, this vanishes
along the edges and has an atom at each vertex.  But sometimes
we wish to view a polygon as an approximation to a smooth curve
and thus spread this curvature out into a smooth density.
For instance the elastic energy, measured as $\int\kappa^2\,ds$ for smooth
curves, blows up if measured directly on a polygon; instead
of squaring~$\K$, we should find a smoothed \emph{curvature density}~$\kappa$
and square that.

For simplicity, we will consider here only the
case of equilateral polygons in~$\E^d$,
where each edge has length~$1$.
To each vertex~$v$, we allocate the neighborhood~$N_v$
consisting of the nearer halves of the
two edges incident to~$v$, with total length~$1$.
Depending on whether we are thinking
of~$\tc$ or~$\tce$, we see total curvature either~$\theta$
or $2\sin(\theta/2)$ at~$v$, and so it would be natural to
use either~$\theta$ or $2\sin(\theta/2)$ as the curvature
density~$\kappa$ along~$N_v$.
The latter has also a geometric interpretation: if $uvw$ are
consecutive vertices of the equilateral polygon, then the circle through
these three points has curvature density $\kappa=2\sin(\theta/2)$.

However, from a number of points of view, there is another
even better measure of curvature density.  Essentially, what
we have ignored above is that when we round off the corners
of a polygon to make a smooth curve, we tend to make the curve
shorter.  Thus, while~$\theta$ or $2\sin(\theta/2)$ might be the
correct total curvature for a neighborhood of~$v$, perhaps it should
get averaged over a length less than~$1$.

Let us consider a particularly simple smoothing, which gives
a $C^{1,1}$ and piecewise circular curve.  Given a polygon~$P$,
replace the neighborhood~$N_v$ of each vertex~$v$
by an ``inscribed'' circular arc, tangent at each endpoint.
This arc turns a total angle~$\theta$, but since it is shorter
than~$N_v$, its curvature density is $\kappa=2\tan(\theta/2)$.

As a simple example, suppose~$P$ is a regular $n$-gon in the plane
with edges of length~$1$ and turning angles $2\pi/n$. 
Its inscribed circle has curvature density $2\tan(\theta/2)$,
while its circumscribed circle
has (smaller) curvature density $2\sin(\theta/2)$.
(Of course, the value~$\theta$ lies between these two.) 

Using the formula $2\tan(\theta/2)$ for the curvature density along
a polygon~$P$ has certain advantages.
For instance if $\kappa(P)\le C$ then we know there is a nearby
$C^{1,1}$ curve (the smoothing by inscribed circular arcs
we used above) with this same curvature bound.
The fact that $2\tan(\theta/2)$ blows up for $\theta=\pi$
reflects the fact that a polygonal corner of turning angle~$\pi$
is really a cusp.  For instance, when the turning angle of a polygon
in the plane passes through~$\pi$, the total signed curvature (or
\emph{turning number}) jumps.  For a smooth curve, such a jump
cannot happen, unless the bending energy blows up because of a cusp.
A bending energy\index{bending energy (discrete)}
for polygons based on $\kappa=2\tan(\theta/2)$
will similarly blow up if we try to change the turning number.

Of course, the bending energy for curves is one conserved quantity
for the for the integrable system related to the Hasimoto\index{Hashimoto flow/surface}
or smoke-ring flow.  In the theory of discrete integrable systems,
it seems clear due to work of Hoffmann and others
that $2\tan(\theta/2)$ is the right notion of curvature density
for equilateral polygons.  See \cite{HoffmannKutz,HoffmannDDG}.

\subsection*{Acknowledgments}
My thoughts on curves of finite total curvature have developed
over the course of writing various collaborative papers
\cite{KS-disto,CKKS,CFKSW1,DS-disto,DDS}
in geometric knot theory.  Thus I owe a great debt to my coauthors
Jason Cantarella, Elizabeth Denne, Yuanan Diao, Joe Fu, Greg Kuperberg,
Rob Kusner and Nancy Wrinkle.
I also gratefully acknowledge helpful conversations with
many other colleagues, including Stephanie Alexander, Dick Bishop,
Mohammad Ghomi and G\"unter M.~Ziegler.

%\bibliographystyle{hamsalpha3}
%\bibliography{thick,drl} 

\begin{thebibliography}{CK{\etalchar{+}}03}

\bibitem[AB98]{AlxBsh:FM}
Stephanie~B. Alexander and Richard~L. Bishop, \emph{The {F}{\'a}ry--{M}ilnor
  theorem in {H}adamard manifolds}, Proc. Amer. Math. Soc. \textbf{126}:11
  (1998), 3427--3436.

\bibitem[AD97]{AyaDub}
Sahbi Ayari and Serge Dubuc, \emph{La formule de {C}auchy sur la longueur d'une
  courbe}, Canad. Math. Bull. \textbf{40}:1 (1997), 3--9.

\bibitem[AFP00]{AFP}
Luigi Ambrosio, Nicola Fusco, and Diego Pallara, \emph{Functions of bounded
  variation and free discontinuity problems}, Clarendon/Oxford, 2000.

\bibitem[AG95]{AltGod}
Helmut Alt and Michael Godau, \emph{Computing the {F}r\'echet distance between
  two polygonal curves}, Internat. J. Comput. Geom. Appl. \textbf{5}:1-2
  (1995), 75--91, Proc. 8th ACM Symp. Comp. Geom. (Berlin, 1992).

\bibitem[AR89]{AlexResh}
Aleksandr~D. Alexandrov and Yuri~G. Reshetnyak, \emph{General theory of
  irregular curves}, Math. Appl. (Soviet) 29, Kluwer, Dordrecht, 1989.

\bibitem[AT72]{AlpTor}
Louis~I. Alpert and Leopoldo~V. Toralballa, \emph{An elementary definition of
  surface area in {$\E^{n+1}$} for smooth surfaces}, Pacific J. Math.
  \textbf{40} (1972), 261--268.

\bibitem[AZ98]{AigZieg}
Martin Aigner and G{\"u}nter~M. Ziegler, \emph{Proofs from {T}he {B}ook},
  Springer, 1998.

\bibitem[Ber98]{Berberian}
Sterling~K. Berberian, \emph{Fundamentals of real analysis}, Springer, 1998.

\bibitem[BH74]{BriHsi}
Frederick Brickell and Chuan-Chih Hsiung, \emph{The total absolute curvature of
  closed curves in {R}iemannian manifolds}, J. Differential Geometry \textbf{9}
  (1974), 177--193.

\bibitem[Bla21]{Blaschke1}
Wilhelm Blaschke, \emph{Ungleichheiten von {H}. {A}. {S}chwarz und {A}. {S}chur
  f{\"u}r {R}aumkurven mit vorgeschriebener {K}r{\"u}mmung.}, Hamb. Abh.
  \textbf{1} (1921), 49--53.

\bibitem[Bla24]{Blaschke2}
\bysame, \emph{Vorlesungen {\"u}ber {D}ifferentialgeometrie},
  Springer, Berlin, 1924.

\bibitem[Boa96]{Boas}
Ralph~P. Boas, Jr., \emph{A primer of real functions}, fourth ed., Carus Math.
  Manus., no.~13, Math. Assoc. Amer., 1996.

\bibitem[Bor47]{Borsuk}
Karol Borsuk, \emph{Sur la courbure totale des courbes ferm\'ees}, Ann. Soc.
  Polon. Math. \textbf{20} (1947), 251--265 (1948).

\bibitem[BS99]{BuckSimon}
Gregory~R. Buck and Jonathon~K. Simon, \emph{Thickness and crossing number of knots}, Topol.
  Appl. \textbf{91}:3 (1999), 245--257.

\bibitem[Buf77]{Buffon}
Georges-Louis~{Leclerc, Comte de} Buffon, \emph{Essai d'arithm{\'e}tique
  morale}, Histoire naturelle, g{\'e}n{\'e}rale er particuli{\`e}re, Suppl.~4,
  1777, pp.~46--123; \url{www.buffon.cnrs.fr/ice/ice_book_detail-fr-text-buffon-buffon_hn-33-7.html}.

\bibitem[Car18]{Carath}
Constantin Carath{\'e}odory, \emph{Vorlesungen {\"u}ber reele {F}unktionen},
  Teubner, 1918, reprinted 2004 by AMS/Chelsea.

\bibitem[Cau13]{Cauchy}
Augustin-Louis Cauchy, \emph{Deuxi{\`e}me m{\'e}moire sur les polygones et les
  poly{\`e}dres}, J. {\'E}cole Polytechnique \textbf{16} (1813), 87--98.

\bibitem[Cau41]{Cauchy-area}
\bysame, \emph{Notes sur divers th{\'e}or{\`e}mes relatifs {\`a}
  la rectification des courbes, et {\`a} la quadrature des surfaces}, C. R.
  Acad. Sci. Paris \textbf{13} (1841), 1060--1063, reprinted in \emph{Oeuvres
  compl{\`e}tes 6}, Gauthier-Villars, 1888, pp.~369--375.

\bibitem[Ces56]{Cesari}
Lamberto Cesari, \emph{Surface area}, Ann. Math. Stud., no.~35, Princeton,
  1956.

\bibitem[Ces89]{Cesari-Chern}
\bysame, \emph{Surface area}, Global Differential Geometry (S.~S.
  Chern, ed.), Math. Assoc. Amer., 1989, pp.~270--302.

\bibitem[CF{\etalchar{+}}04]{CFKSW1}
Jason Cantarella, Joseph~H.~G. Fu, Robert~B. Kusner,
 John~M. Sullivan, and Nancy~C. Wrinkle,
  \emph{Criticality for the {G}ehring link problem}, Geometry and Topology
  \textbf{10} (2006), 2055--2115; \url{arXiv.org/math.DG/0402212}. 

\bibitem[Cha62]{Chak1}
Gulbank~D. Chakerian, \emph{An inequality for closed space curves}, Pacific J.
  Math. \textbf{12} (1962), 53--57.

\bibitem[Cha64]{Chak2}
\bysame, \emph{On some geometric inequalities}, Proc. Amer. Math.
  Soc. \textbf{15} (1964), 886--888.

\bibitem[Che89]{Chern}
Shiing~Shen Chern, \emph{Curves and surfaces in euclidean space}, Global
  Differential Geometry (S.~S. Chern, ed.), Math. Assoc. Amer., 1989,
  pp.~99--139.

\bibitem[CK{\etalchar{+}}03]{CKKS}
Jason Cantarella, Greg Kuperberg, Robert~B. Kusner, and John~M. Sullivan,
  \emph{The second hull of a knotted curve}, Amer. J. Math \textbf{125} (2003),
  1335--1348; \url{arXiv.org/math.GT/0204106}.

\bibitem[CKS02]{CKS2}
Jason Cantarella, Robert~B. Kusner, and John~M. Sullivan, \emph{On the minimum
  ropelength of knots and links}, Invent. Math. \textbf{150} (2002), 257--286;
  \url{arXiv.org/math.GT/0103224}.

\bibitem[Con82]{Connelly-invent}
Robert Connelly, \emph{Rigidity and energy}, Invent. Math. \textbf{66} (1982),
  11--33.

\bibitem[Cro68]{Crofton}
Morgan~W. Crofton, \emph{On the theory of local probability}, Phil. Trans. R.
  Soc. London \textbf{158} (1868), 181--199.

\bibitem[Cro97]{Cromwell-polyh}
Peter Cromwell, \emph{Polyhedra}, Cambridge, 1997.

\bibitem[dC76]{doCarmo}
Manfredo~P. do~Carmo, \emph{Differential geometry of curves and surfaces},
  Prentice-Hall, 1976.

\bibitem[DDS06]{DDS}
Elizabeth Denne, Yuanan Diao, and John~M. Sullivan, \emph{Quadrisecants give
  new lower bounds for the ropelength of a knot}, Geometry and Topology
  \textbf{10} (2006), 1--26; \url{arXiv.org/math.DG/0408026}.

\bibitem[Den04]{Denne}
Elizabeth Denne, \emph{Alternating quadrisecants of knots}, Ph.D. thesis, Univ.
  Illinois, Urbana, 2004, \url{arXiv.org/math.GT/0510561}.

\bibitem[DS04]{DS-disto}
Elizabeth Denne and John~M. Sullivan, \emph{The distortion of a knotted curve},
  preprint, 2004, \url{arXiv.org/math.GT/0409438v1}.

\bibitem[DS08]{DS-ConvFTC}
\bysame,
\emph{Convergence and isotopy type for graphs of finite total curvature},
\thisvol{Denne}; \url{arXiv.org/math.GT/0606008}.

\bibitem[F{\'a}r49]{fary}
Istv{\'a}n F{\'a}ry, \emph{Sur la courbure totale d'une courbe gauche faisant
  un n{\oe}ud}, Bull. Soc. Math. France \textbf{77} (1949), 128--138.

\bibitem[Fen29]{Fenchel}
Werner Fenchel, \emph{{\"U}ber {K}r{\"u}mmung und {W}indung geschlossener
  {R}aumkurven}, Math. Ann. \textbf{101} (1929), 238--252;
  \url{www-gdz.sub.uni-goettingen.de/cgi-bin/digbib.cgi?PPN235181684_0101}.

\bibitem[Fen51]{MR12:634d}
Werner Fenchel, \emph{On the differential geometry of closed space curves},
  Bull. Amer. Math. Soc. \textbf{57} (1951), 44--54.

\bibitem[Fr{\'e}05]{Frechet}
Maurice Fr{\'e}chet, \emph{Sur l'{\'e}cart de deux courbes et sur les courbes
  limites}, Trans. Amer. Math. Soc. \textbf{6} (1905), 435--449.

\bibitem[GP83]{GoffPed}
Casper Goffman and George Pedrick, \emph{First course in functional analysis},
  second ed., Chelsea, 1983.

\bibitem[Gra46]{Graves}
Lawrence~M. Graves, \emph{Theory of functions of real variables}, McGraw Hill,
  1946.

\bibitem[Gro78]{GromovHED}
Mikhael Gromov, \emph{Homotopical effects of dilatation}, J. Diff. Geom.
  \textbf{13} (1978), 303--310.

\bibitem[Gro81]{GLP}
\bysame, \emph{Structures m\'etriques pour les vari\'et\'es
  riemanniennes}, Cedic, Paris, 1981, edited by J. Lafontaine and P. Pansu.

\bibitem[Gro83]{GromovFRM}
\bysame, \emph{Filling {R}iemannian manifolds}, J.~Diff. Geom.
  \textbf{18} (1983), 1--147.

\bibitem[HK04]{HoffmannKutz}
Tim Hoffmann and Nadja Kutz, \emph{Discrete curves in
  {$\mathbb{C}\mathrm{P}^1$} and the {T}oda lattice}, Stud. Appl. Math.
  \textbf{113}:1 (2004), 31--55; \url{arXiv.org/math.DG/0208190}.

\bibitem[Hof08]{HoffmannDDG}
Tim Hoffmann,
\emph{Discrete Hashimoto surfaces and a doubly discrete smoke-ring flow},
\thisvol{Hoffmann}; \url{arXiv.org/math.DG/0007150}.

\bibitem[Hor71]{Horn-Fenchel}
Roger~A. Horn, \emph{On {F}enchel's theorem}, Amer. Math. Monthly \textbf{78}
  (1971), 380--381.

\bibitem[Jor93]{Jordan}
Camille Jordan, \emph{Cours d'analyse de l'{\'e}cole polytechnique},
  Gauthier-Villars, 1893.

\bibitem[KR97]{KlainRota}
Daniel~A. Klain and Gian-Carlo Rota, \emph{Introduction to geometric
  probability}, Cambridge, 1997.

\bibitem[KS97]{KS-disto}
Robert~B. Kusner and John~M. Sullivan, \emph{On distortion and thickness of
  knots}, Topology and Geometry in Polymer Science (Whittington, Sumners, and
  Lodge, eds.), IMA Vol.~103, Springer, 1997, pp.~67--78;
  \url{arXiv.org/dg-ga/9702001}.

\bibitem[Leb02]{Lebesgue}
Henri Lebesgue, \emph{Int\'egrale, longuer, aire}, Annali di Mat. pura appl.
  \textbf{7} (1902), 231--359.

\bibitem[Lie29]{Liebmann-Fenchel}
Heinrich Liebmann, \emph{Elementarer {B}eweis des {F}enchelschen {S}atzes
  {\"u}ber die {K}r{\"u}mmung geschlossener {R}aumkurven}, Sitz.ber. Akad.
  Berlin (1929), 392--393.

\bibitem[Mat03]{Matousek}
Ji{\v{r}}{\'\i} Matou{\v{s}}ek, \emph{Using the {B}orsuk--{U}lam theorem},
  Springer, Berlin, 2003.

\bibitem[Mil50]{Milnor}
John~W. Milnor, \emph{On the total curvature of knots}, Ann. of Math.
  \textbf{52} (1950), 248--257.

\bibitem[Mor36]{Morse-param}
Marston Morse, \emph{A special parametrization of curves}, Bull. Amer. Math.
  Soc. \textbf{42} (1936), 915--922.

\bibitem[Mor88]{morgan}
Frank Morgan, \emph{Geometric measure theory: A beginner's guide}, Academic
  Press, 1988.

\bibitem[O'R00]{ORourke-Schur}
Joseph O'Rourke,
 \emph{On the development of the intersection of a plane with a polytope},
 Comput. Geom. Theory Appl. \textbf{24}:1 (2003), 3--10;
 \url{arXiv.org/cs.CG/0006035v3}.

\bibitem[San89]{Santalo-Chern}
Luis~A. Santal{\'o}, \emph{Integral geometry}, Global Differential Geometry,
  Math. Assoc. Amer., 1989, pp.~303--350.

\bibitem[San04]{Santalo}
\bysame, \emph{Integral geometry and geometric probability}, second
  ed., Cambridge, 2004.

\bibitem[Sch85]{Scheeffer}
Ludwig Scheeffer, \emph{Allgemeine {U}ntersuchungen {\"u}ber {R}ectification
  der {C}urven}, Acta Math. \textbf{5} (1884--85), 49--82.

\bibitem[Sch90]{Schwarz-lantern}
Hermann~Amandus Schwarz, \emph{Sur une d{\'e}finition erron{\'e}e de l'aire
  d'une surface courbe}, Ges. math. Abhandl., vol.~2, Springer, 1890,
  pp.~309--311 and 369--370.

\bibitem[Sch21]{Schur}
Axel Schur, \emph{{\"U}ber die {S}chwarzsche {E}xtremaleigenschaft des
  {K}reises unter den {K}urven konstanter {K}r{\"u}mmung}, Math. Annalen
  \textbf{83} (1921), 143--148;
  \url{www-gdz.sub.uni-goettingen.de/cgi-bin/digbib.cgi?PPN235181684_0083}.

\bibitem[Sch25]{Schmidt-Schur}
Erhard Schmidt, \emph{{\"U}ber das {E}xtremum der {B}ogenl{\"a}nge einer
  {R}aumkurve bei vor\-ge\-schriebenen {E}inschr{\"a}nkungen ihrer {K}r{\"u}mmung},
  Sitz.ber. Akad. Berlin (1925), 485--490.

\bibitem[Sch98]{Schmitz}
Carsten Schmitz, \emph{The theorem of {F}{\'a}ry and {M}ilnor for {H}adamard
  manifolds}, Geom. Dedicata \textbf{71}:1 (1998), 83--90.

\bibitem[Ser68]{Serret}
Joseph-Alfred Serret, \emph{Cours de calcul diff{\'e}rentiel et int{\'e}gral},
  vol.~2, Gauthier-Villars,\linebreak Paris, 1868;\hspace{1em}
  \url{historical.library.cornell.edu/cgi-bin/cul.math/docviewer?did=05270001&seq=311}.

\bibitem[Sul08]{Sul-CSDS}
John~M.~Sullivan, \emph{Curvatures of smooth and discrete surfaces},
\thisvol{SullivanCSDS}; \url{arXiv.org/0710.4497}.

\bibitem[Tan98]{Taniyama-tcgraph}
Kouki Taniyama, \emph{Total curvature of graphs in Euclidean spaces},
  Differential Geom. Appl. \textbf{8}:2 (1998), 135--155.

\bibitem[Ton21]{Tonelli}
Leonida Tonelli, \emph{Fondamenti di calcolo delle variazione}, Zanichelli,
  Bologna, 1921.

\bibitem[Vos55]{Voss-Fenchel}
Konrad Voss, \emph{Eine {B}emerkung \"uber die {T}otalkr\"ummung geschlossener
  {R}aumkurven}, Arch. Math. \textbf{6} (1955), 259--263.

\bibitem[War08]{xWardetzky}
Max Wardetzky, \emph{Convergence of the cotangent formula: An overview},
\thisvol{Wardetzky}.

\bibitem[Zie89]{Ziemer}
William~P. Ziemer, \emph{Weakly differentiable functions}, GTM 120, Springer,
  1989.

\end{thebibliography}
\newcommand{\etalchar}[1]{$^{#1}$}
\providecommand{\bysame}{\leavevmode\hbox to3em{\hrulefill}\thinspace}

%\end{document}
\end{document}